\documentclass[10pt,reqno]{amsart}
\usepackage{mathrsfs}
\usepackage{amssymb}
\usepackage{amsmath}
\usepackage{amsfonts}
\usepackage{amscd}
\usepackage{graphicx}
\usepackage[colorlinks, linkcolor=blue, anchorcolor=blue, citecolor=blue]{hyperref}
\begin{document}
\setlength{\oddsidemargin}{0cm} \setlength{\evensidemargin}{0cm}
\baselineskip=20pt

\theoremstyle{plain} \makeatletter
\newtheorem{theorem}{Theorem}[section]
\newtheorem{proposition}[theorem]{Proposition}
\newtheorem{lemma}[theorem]{Lemma}
\newtheorem{coro}[theorem]{Corollary}

\theoremstyle{definition}
\newtheorem{defi}[theorem]{Definition}
\newtheorem{notation}[theorem]{Notation}
\newtheorem{exam}[theorem]{Example}
\newtheorem{prop}[theorem]{Proposition}
\newtheorem{conj}[theorem]{Conjecture}
\newtheorem{prob}[theorem]{Problem}
\newtheorem{remark}[theorem]{Remark}
\newtheorem{claim}{Claim}

\newcommand{\SO}{{\mathrm S}{\mathrm O}}
\newcommand{\SU}{{\mathrm S}{\mathrm U}}
\newcommand{\Sp}{{\mathrm S}{\mathrm p}}
\newcommand{\s}{{\mathfrak s}{\mathfrak p}}
\newcommand{\so}{{\mathfrak s}{\mathfrak o}}
\newcommand{\Ad}{{\mathrm A}{\mathrm d}}
\newcommand{\ad}{{\mathrm a}{\mathrm d}}
\newcommand{\id}{{\mathrm i}{\mathrm d}}
\newcommand{\E}{{\mathrm E}}
\newcommand{\F}{{\mathrm F}}
\newcommand{\G}{{\mathrm G}}
\newcommand{\R}{{\mathrm Ric}}

\numberwithin{equation}{section}

\title{Notes on ``Einstein metrics on compact simple Lie groups attached to standard triples"}
\author{Huibin Chen}
\address{School of Mathematical Sciences and LPMC, Nankai University,
Tianjin 300071, P.R. China}
\author{Zhiqi Chen}
\address{School of Mathematical Sciences and LPMC, Nankai University,
Tianjin 300071, P.R. China}\email{chenzhiqi@nankai.edu.cn}
\date{}
\maketitle
\begin{abstract}
In the paper ``Einstein metrics on compact simple Lie groups attached to standard triples", the authors introduced the definition of standard triples and proved that every compact simple Lie group $G$ attached to a standard triple $(G,K,H)$ admits a left-invariant Einstein metric which is not naturally reductive except the standard triple $(\Sp(4),2\Sp(2),4\Sp(1))$. For the triple $(\Sp(4),2\Sp(2),4\Sp(1))$, we find there exists an involution pair of $\s(4)$ such that $4\s(1)$ is the fixed point of the pair, and then give the decomposition of $\s(4)$ as a direct sum of irreducible $\ad(4\s(1))$-modules. But $\Sp(4)/4\Sp(1)$ is not a generalized Wallach space. Furthermore we give left-invariant Einstein metrics on $\Sp(4)$ which are non-naturally reductive and $\Ad(4\Sp(1))$-invariant. For the general case $(\Sp(2n_1n_2),2\Sp(n_1n_2),2n_2\Sp(n_1))$, there exist $2n_2-1$ involutions of $\s(2n_1n_2)$ such that $2n_2\s(n_1))$ is the fixed point of these $2n_2-1$ involutions, and it follows the decomposition of $\s(2n_1n_2)$ as a direct sum of irreducible $\ad(2n_2\s(n_1))$-modules. In order to give new non-naturally reductive and $\Ad(2n_2\Sp(n_1)))$-invariant Einstein metrics on $\Sp(2n_1n_2)$, we prove a general result, i.e. $\Sp(2k+l)$ admits at least two non-naturally reductive Einstein metrics which are $\Ad(\Sp(k)\times\Sp(k)\times\Sp(l))$-invariant if $k<l$. It implies that every compact simple Lie group $\Sp(n)$ for $n\geq 4$ admits at least $2[\frac{n-1}{3}]$ non-naturally reductive left-invariant Einstein metrics.
\end{abstract}

\section{Introduction}
A Riemannian manifold $(M,g)$ is called Einstein if there exists a constant $\lambda\in\mathbb{R}$ such that the Ricci tensor $r$ with respect to $g$ satisfies $r=\lambda g$. See Besse's book \cite{Be}, the papers of Jensen \cite{GJ1973} and Wang and Ziller \cite{WaZi} for more details. For Lie groups, D'Atri and Ziller prove in \cite{JW} that every compact simple Lie group except $\SO(3)$ admits at least two left-invariant Einstein metrics which are naturally reductive. But they also mention in \cite{JW} that it is difficult to give non-naturally reductive Einstein metrics on compact simple Lie groups.

Recently, there are some studies on non-naturally reductive Einstein metrics on compact simple Lie groups. Mori gives in \cite{Mo} a class of non-naturally reductive left-invariant Einstein metrics on compact simple Lie groups $\SU(n)$ for $n\geq 6$. After that, Arvanitoyeorgos, Mori and Sakane prove the existence of non-naturally reductive Einstein metrics on $\SO(n)(n\geq11)$, $\Sp(n)(n\geq3)$, $\E_6$, $\E_7$ and $\E_8$. Then Chen and Liang obtain in \cite{ChLi} a non-naturally reductive Einstein metric on the compact simple Lie group $\F_4$. Recently, Chrysikos and Sakane find new non-naturally reductive Einstein metrics on exceptional Lie groups in \cite{ChSa}, especially they give the first non-naturally reductive Einstein metric on $\G_2$.

In \cite{YD1}, Yan and Deng study non-naturally reductive Einstein metrics on a compact Lie group $G$ associated with a triple $(G,K,H)$. Here $G$ is a compact simple Lie group, $K$ is a closed subgroup of $G$ such that $G/K$ is a compact irreducible symmetric space, and $H$ is a closed subgroup of $K$. Denote the Lie algebras of $G,K,H$ by $\mathfrak g,\mathfrak k,\mathfrak h$ respectively. Let $B,B_{\mathfrak k},B_{\mathfrak h}$ be the negative of the Killing forms of $\mathfrak g,\mathfrak k,\mathfrak h$ respectively. Consider the $B$-orthogonal decomposition
\begin{equation}\label{decomp}\mathfrak g=\mathfrak k\oplus\mathfrak p=\mathfrak h\oplus\mathfrak u\oplus\mathfrak p \end{equation} and the left-invariant metrics on $G$ determined by the $\ad({\mathfrak h})$-invariant metric on $\mathfrak g$ of the form \begin{equation}\label{metric-form}\langle\cdot,\cdot\rangle=B|_{\mathfrak h}+xB|_{\mathfrak u}+yB|_{\mathfrak p},\quad x,y\in\mathbb R^+.\end{equation}
A triple $(G,K,H)$ (or ($\mathfrak g$,$\mathfrak k$,$\mathfrak h$)) is called a basic triple in \cite{YD1} if
\begin{enumerate}
   \item $B_{\mathfrak k}=c_1B|_{\mathfrak k}$, $B_{\mathfrak h}=c_2B|_{\mathfrak h}$ for some $c_1>0,c_2>0$, and
   \item $\sum_i(\ad^2(h_i))|_{\mathfrak u}=-\lambda_u\id$, $\sum_i(\ad^2(h_i))|_{\mathfrak p}=-\lambda_p\id$ for some $\lambda_u>0,\lambda_p>0$. Here $\{h_i\}$ is a $B$-orthonormal basis of $\mathfrak h$.
\end{enumerate}
A basic triple $(G,K,H)$ is called standard if the standard homogenous metric $g_B$ on $G/H$ is Einstein. In fact, the classification of the homogeneous space whose standard homogeneous metric is Einstein is given by Wang and Ziller in \cite{WZ1} . Thus it is easy to give a complete classification of standard triples based on \cite{WZ1,Wo1}. Furthermore in \cite{YD1}, Yan and Deng classify standard triples and prove that every compact simple Lie group $G$ attached to a standard triple $(G,K,H)$ admits a left-invariant Einstein metric of the form~\eqref{metric-form} which is non-naturally reductive except the standard triple $(\Sp(4),2\Sp(2),4\Sp(1))$.

In general, for a basic (even standard) triple $(G,K,H)$, the decomposition~\eqref{decomp} is not necessary to be a direct sum of irreducible $\ad({\mathfrak h})$-modules. That is, the metric of the form~\eqref{metric-form} is only a special class of left-invariant metrics on $G$ which are also $\Ad(H)$-invariant. In order to obtain non-naturally reductive left-invariant Einstein metrics on $\Sp(4)$ which are $\Ad(4\Sp(1))$-invariant, we need to discuss the structure of $\s(4)$ as a direct sum of irreducible $\ad({4\s(1)})$-modules. Firstly, we have the following theorem.
\begin{theorem}\label{structure}
For the standard triple $(\Sp(4),2\Sp(2),4\Sp(1))$, there exists an involution pair $(\theta,\tau)$ of $\s(4)$ such that $\theta\tau=\tau\theta$ and $4\s(1)=\{x\in \s(4)|\theta(x)=x,\tau(x)=x\}.$
\end{theorem}

It follows the decomposition of $\s(4)$ as a direct sum of irreducible $\ad({4\s(1)})$-modules. Although $4\s(1)$ is the fixed point of an involution pair of $\s(4)$, $\Sp(4)/4\Sp(1)$ is not a generalized Wallach space, which is also called a three-locally-symmetric space. The notation of a generalized Wallach space is introduced by Nikonorov in \cite{Ni00}. A compact homogeneous space $G/H$ with a semisimple connected Lie group $G$ and a connected Lie subgroup $H$ is a generalized Wallach space if $\mathfrak m$ is the direct sum of three irreducible $\ad{\mathfrak h}$-modules pairwise orthogonal with respect to $B$, i.e.
$$\mathfrak m=\mathfrak m_1\oplus \mathfrak m_2\oplus\mathfrak m_3,$$
with $[\mathfrak m_i,\mathfrak m_i]\subset \mathfrak h$ for any $i\in \{1,2,3\}$. Here $B$ is the negative of Killing from of $\mathfrak g$, $\mathfrak h$ is the Lie algebra of $H$, and $\mathfrak m$ is the orthogonal complement of $\mathfrak h$ in $\mathfrak g$ with respect to $B$. The classification of generalized Wallach spaces is given by Nikonorov in \cite{Ni16}. In \cite{ChKaLi}, Chen, Kang and Liang prove that the classification of generalized Wallach spaces is equivalent with the classification of involution pairs of compact Lie groups satisfying certain conditions, and then give the classification for compact simple Lie groups. 

In \cite{ArSaSt}, Arvanitoyeorgos, Sakane and Statha study non-naturally reductive Einstein metrics on $\Sp(k_1+k_2+k_3)$ which are $\Ad(\Sp(k_1)\times\Sp(k_2)\times\Sp(k_3))$-invariant. In particular, they prove that $\Sp(k_1+2)$ admits a non-naturally reductive Einstein metric which is  
$\Ad(\Sp(k_1)\times\Sp(1)\times\Sp(1))$-invariant, and give three non-naturally reductive Einstein metrics for $k_1=2$. Based on the decomposition of $\s(4)$ as a direct sum of irreducible $\ad({4\s(1)})$-modules and the results in \cite{ArSaSt}, we have the following theorem.

\begin{theorem}\label{admitmetric}
The compact simple Lie group $\Sp(4)$ admits non-naturally reductive left-invariant Einstein metrics which are $\Ad(4\Sp(1))$-invariant.
\end{theorem}

For the standard triple $(\Sp(2n_1n_2),2\Sp(n_1n_2),2n_2\Sp(n_1))$ except the case $(\Sp(4),2\Sp(2),4\Sp(1))$, Yan and Deng prove in \cite{YD1} that $\Sp(2n_1n_2)$ attached to the triple $(\Sp(2n_1n_2),2\Sp(n_1n_2),2n_2\Sp(n_1))$ admits a non-naturally reductive left-invariant Einstein metric of the form~\eqref{metric-form} which is $\Ad(2n_2\Sp(n_1))$-invariant. Furthermore, they obtain an interesting result on the number of Einstein metrics.

\begin{lemma}[\cite{YD1}]\label{lemma}
 For any integer $n=p_1^{l_1}p_2^{l_2}\cdots p_s^{l_s}$ with $p_i$ prime and $p_i\not=p_j$, $\Sp(2n)$ admits at least $$(l_1+1)(l_2+1)\cdots (l_s+1)-1$$ non-equivalent non-naturally reductive Einstein metrics.
\end{lemma}

In this paper, we character the structure of $\s(2n_1n_2)$ as a direct sum of irreducible $\ad(2n_2\s(n_1)))$-modules by the following theorem.
\begin{theorem}\label{structure1}
For the standard triple $(\Sp(2n_1n_2),2\Sp(n_1n_2),2n_2\Sp(n_1))$, there exists involutions $\theta_i,1\leq i\leq 2n_2-1$ of $\s(2n_1n_2)$ such that $\theta_i\theta_j=\theta_j\theta_i$ and $2n_2\s(n_1)=\{x\in \s(2n_1n_2)|\theta_i(x)=x,1\leq i\leq 2n_2-1\}.$
\end{theorem}

It follows from Theorem~\ref{structure1} that a $B$-orthogonal decomposition of $\s(2n_1n_2)$ as a direct sum of irreducible $\ad(2n_2\s(n_1))$-modules. Here we point out that, from Theorem~\ref{structure1}, we get another decomposition of $\s(4)$ by three involutions of $\s(4)$ directly, which is in essential the same as that from Theorem~\ref{structure}. 

In order to give new non-naturally reductive Einstein metrics on $\Sp(2n_1n_2)$ which are $\Ad(2n_2\Sp(n_1))$-invariant, we consider another class of left-invariant metrics on $\Sp(2n_1n_2)$ which are $\Ad(2n_2\Sp(n_1))$-invariant different from those given in \cite{YD1}. The metrics correspond to the decomposition of $\Sp(2n_1n_2)$ under two special involutions from $\theta_i,1\leq i\leq 2n_2-1$ whose fixed point is $\s(n_1)\oplus\s(n_1)\oplus \s(2(n_2-1)n_1)$. More general, we prove the following theorem.

\begin{theorem}\label{theorem}
For given positive integers $k$ and $l$ with $k<l$, $\Sp(2k+l)$ admits at least two non-naturally reductive left-invariant Einstein metrics which are $\Ad(\Sp(k)\times\Sp(k)\times\Sp(l))$-invariant.
\end{theorem}

In particular, in Theorem~\ref{theorem}, let $n_1=k$ and take $l=2(n_2-1)k$, we have at least two non-naturally reductive left-invariant Einstein metrics on $\Sp(2n_1n_2)$ which are $\Ad(\Sp(n_1)\times\Sp(n_1)\times\Sp(2(n_2-1)n_1))$. In particular, they are $\Ad(2n_2\Sp(n_1))$-invariant.

Comparing to Lemma~\ref{lemma}, we get the lower bound of the number of non-naturally Einstein metrics on $\Sp(n)$ for $n\geq 4$ from Theorem~\ref{theorem}.

\begin{theorem}\label{maintheorem}
For every $n\geq 4$, the compact simple Lie group $\Sp(n)$ admits at least $2[\frac{n-1}{3}]$ non-naturally reductive left-invariant Einstein metrics.
\end{theorem}

The paper is organized as follows. In section 2, we list some preliminaries which is necessary for this paper. In section 3, we recall the results on the theory of involutions on compact simple Lie groups, and then decompose $\s(4)$ as a direct sum of irreducible $\ad({4\s(1)})$-modules based on the above theory. That is, we prove theorem~\ref{structure}. Then in section 4, we prove Theorem~\ref{admitmetric} by the studies on $\Sp(n)$ given in~\cite{ArSaSt}. In section 5, we decompose $\s(2n_1n_2)$ as a direct sum of irreducible $\ad(2n_2\s(n_1)))$-modules corresponding to $2n_2-1$ involutions of $\s(2n_1n_2)$. In particular, we prove Theorem~\ref{structure1}. Furthermore, we prove Theorems~\ref{theorem} and~\ref{maintheorem} in section 6.

\section{Preliminaries: naturally reductive metrics and the Ricci tensor}
Let $G$ be a compact simple Lie group and let $K$ be a connected closed subgroup of $G$ with Lie algebras $\mathfrak g$ and $\mathfrak k$ respectively. Let $\mathfrak g=\mathfrak k\oplus\mathfrak m$ be the $B$-orthogonal decomposition. Here $[\mathfrak k,\mathfrak m]\subset\mathfrak m$. Assume that $\mathfrak m$ can be decomposed into mutually non-equivalent irreducible $\Ad(K)$-modules:
\begin{equation*}
\mathfrak m=\mathfrak m_1\oplus\cdots\oplus\mathfrak m_q.
\end{equation*}
Let $\mathfrak k=\mathfrak k_0\oplus\mathfrak k_1\oplus\cdots\oplus\mathfrak k_p$, where $\mathfrak k_0=Z(\mathfrak k)$ is the center of $\mathfrak k$ and every $\mathfrak k_i$ is simple for $i=1,\cdots,p$. It is well-known that there exists a one-to-one corresponding between $G$-invariant metrics on $G/K$ and $\Ad(K)$-invariant inner products on $\mathfrak m$. The $G$-invariant metric $\langle\cdot,\cdot\rangle$ on $M=G/K$ is called naturally reductive if
\begin{equation*}
\langle [X,Y]_{\mathfrak m},Z\rangle+\langle Y,[X,Z]_{\mathfrak m}\rangle=0, \quad \forall X,Y,Z\in\mathfrak m.
\end{equation*}

In \cite{JW}, D'Atri and Ziller give a sufficient and necessary condition for a left-invariant metric on a compact simple Lie group to be naturally reductive.
\begin{lemma}[\cite{JW}]\label{nr}
For any inner product $b$ on $\mathfrak k_0$, every left-invariant metric on $G$ with the form
\begin{equation*}
\langle\cdot,\cdot \rangle=u_0b|_{\mathfrak k_0}+u_1B|_{\mathfrak k_1}+\cdots+u_pB|_{\mathfrak k_p}+xB|_{\mathfrak m},\quad(u_0,u_1,\cdots,u_p,x\in\mathbb{R}^+)
\end{equation*}
is naturally reductive with respect to the action $(g,k)y=gyk^{-1}$ of $G\times K$.
Conversely, if a left-invariant metric $\langle\cdot,\cdot\rangle$ on a compact simple Lie group $G$ is naturally reductive, then there exists a closed subgroup $K$ of $G$ such that $\langle\cdot,\cdot\rangle$ can be written as above.
\end{lemma}

Now we have a $B$-orthogonal decomposition of $\mathfrak g$ which is $\Ad(K)$-invariant: $$\mathfrak g=\mathfrak k_0\oplus\mathfrak k_1\oplus\cdots\oplus\mathfrak k_p\oplus\mathfrak m_1\oplus\cdots\oplus\mathfrak m_q=\mathfrak k_0\oplus\mathfrak k_1\oplus\cdots\oplus\mathfrak k_p\oplus\mathfrak k_{p+1}\oplus\cdots\oplus\mathfrak k_{p+q},$$ where $\mathfrak k_{p+i}=\mathfrak m_i$ for $1\leq i\leq q$. Assume that $\mathrm{dim}\mathfrak k_0\leq 1$. Consider the following left-invariant metric on $G$ which is $\Ad(K)$-invariant:
\begin{equation}\label{metric1}
\langle\cdot,\cdot\rangle=x_0\cdot B|_{\mathfrak k_0}+x_1\cdot B|_{\mathfrak k_1}+\cdots+x_{p+q}\cdot B|_{\mathfrak k_{p+q}},
\end{equation}
where $x_i\in\mathbb{R}^+$ for $i=1,\cdots,p+q$.

Let $d_i=\dim\mathfrak k_i$ and let $\{e^{i}_{\alpha}\}_{\alpha=1}^{d_i}$ be a $B$-orthonormal basis of $\mathfrak k_i$. Denote  $A_{\alpha,\beta}^{\gamma}=B([e_{\alpha}^{i},e_{\beta}^{j}],e_{\gamma}^{k})$, i.e.  $[e_{\alpha}^{i},e_{\beta}^{j}]=\sum_{\gamma}A_{\alpha,\beta}^{\gamma}e_{\gamma}^{k}$, and define
\begin{equation*}
(ijk):=\sum(A_{\alpha,\beta}^{\gamma})^2,
\end{equation*}
where the sum is taken over all indices $\alpha,\beta,\gamma$ with $e_{\alpha}^{i}\in\mathfrak k_i, e_{\beta}^{j}\in\mathfrak k_j, e_{\gamma}^{k}\in\mathfrak k_k$. Then $(ijk)$ is independent of the choice for the $B$-orthonormal basis of $\mathfrak k_i,\mathfrak k_j,\mathfrak k_k$, and $(ijk)=(jik)=(jki)$. Furthermore in \cite{ArMoSa}, there are fundamental formulae for the Ricci tensor for compact Lie groups and compact homogeneous spaces. Here we just give the formula for Lie groups.

\begin{lemma}\cite{ArMoSa,PaSa}\label{formula1}
Let $G$ be a compact connected simple Lie group endowed with a left-invariant metric $\langle\cdot ,\cdot\rangle$ of the form (\ref{metric1}). Then the components $r_0,r_1,\cdots,r_{p+q}$ of the Ricci tensor associated to $\langle\cdot,\cdot\rangle$ are:
\begin{equation*}
r_k=\frac{1}{2x_k}+\frac{1}{4d_k}\sum_{j,i}\frac{x_k}{x_jx_i}(kji)-\frac{1}{2d_k}\sum_{j,i}\frac{x_j}{x_kx_i}(jki),\quad(k=0,1,\cdots,p+q).
\end{equation*}
Here, the sums are taken over all $i=0,1,\cdots,p+q$. In particular for each $k$, $\sum_{i,j}(jki)=d_k.$
\end{lemma}

\section{The decomposition of $\s(4)$ as a direct sum of irreducible $\ad({4\s(1)})$-modules.}\label{sec2}

This section is to character the decomposition of $\s(4)$ as a direct sum of irreducible $\ad({4\s(1)})$-modules according to the standard triple $(\Sp(4),2\Sp(2),4\Sp(1))$. Here, $\Sp(4)/2\Sp(2)$ is a compact irreducible symmetric space. That is, there is an involution $\theta$ of $\s(4)$ such that $$2\s(2)=\{x\in \s(4)|\theta(x)=x\}, \quad\mathfrak p=\{x\in \s(4)|\theta(x)=-x\}, \quad\s(4)=2\s(2)\oplus\mathfrak p.$$

In general, let $G$ be a compact simple connected Lie group with the Lie algebra
${\mathfrak g}$ and let $\theta$ be an involution of $G$. Then for the involuation $\theta$, we have a
decomposition,
$${\mathfrak g}={\mathfrak k}\oplus{\mathfrak p},$$ where ${\mathfrak
k}=\{X\in {\mathfrak g}|\theta(X)=X\}$ and ${\mathfrak p}=\{X\in
{\mathfrak g}|\theta(X)=-X\}$. Cartan and Gantmacher made great attributions on the classification of involutions on compact Lie groups. Let ${\mathfrak t_1}$ be a Cartan subalgebra of ${\mathfrak k}$ and let $\mathfrak t$ be a Cartan subalgebra of
${\mathfrak g}$ containing $\mathfrak t_1$.
\begin{lemma}[Gantmacher Theorem]\label{chooseH} With the above notations, $\theta$ is conjugate with $\theta_0e^{\mathrm{ad}H}$ under
$\mathrm{Aut} {\mathfrak g}$, where $H\in {\mathfrak t_1}$ and $\theta_0$ is an involution which keeps the Dynkin diagram invariant.
\end{lemma}
Let $\Pi=\{\alpha_1,\dots,\alpha_n\}$ be a fundamental system of
${\mathfrak t}$ and $\phi=\sum_{i=1}^nm_i\alpha_i$ be the maximal
root respectively. Let $\alpha_i'=\frac{1}{2}(\alpha_i+\theta_0(\alpha_i))$. Then $\Pi'=\{\alpha_1',\dots,\alpha_l'\}$ consisting different elements in $\{\alpha_1',\dots,\alpha_n'\}$ is a fundamental system of $\mathfrak g_0$, where $\mathfrak g_0=\{X\in {\mathfrak g}|\theta_0(X)=X\}$. Denote by $\phi'=\sum_{i=1}^l m_i'\alpha_i'$ the maximal root of $\mathfrak g_0$ respectively. Furthermore we have
\begin{lemma}[\cite{Ya1}]\label{choosealpha} If $H\not=0$, then for some $i$, we can take $H$ satisfying
\begin{equation}\label{H}
\alpha_i'=\alpha_i; \quad \langle H,\alpha_i'\rangle=\pi\sqrt{-1};\quad \langle
H,\alpha_j'\rangle=0, \forall j\not=i.\end{equation} Here $m_i'=1$ or
$m_i'=2$.
\end{lemma}

Moreover, $\mathfrak k$ is described as follows.

\begin{lemma}[\cite{Ya1}]\label{subalgebra}
Let the notations be as above. Assume that $\alpha_i$ satisfies the
identity~(\ref{H}).
\begin{enumerate}
  \item If $\theta_0=Id$ and $m_i=1$, then $\Pi-\{\alpha_i\}$ is a fundamental
  system of $\mathfrak k$, and $\phi$ and $-\alpha_i$ are the
  highest weights of $\mathrm{ad}_{\mathfrak m^{\mathbb C}}{\mathfrak k}$ corresponding to the fundamental system.
  \item If $\theta_0=Id$ and $m_i=2$, then $\Pi-\{\alpha_i\}\cup\{-\phi\}$ is a fundamental
  system of $\mathfrak k$, and $-\alpha_i$ is the
  highest weight of $\mathrm{ad}_{\mathfrak m^{\mathbb C}}{\mathfrak k}$ corresponding to the fundamental system.
  \item If $\theta_0\not=Id$, then $\Pi'-\{\alpha_i'\}\cup\{\beta_0\}$ is a fundamental
  system of $\mathfrak k$, and $-\alpha_i$ is the
  highest weight of $\mathrm{ad}_{\mathfrak m^{\mathbb C}}{\mathfrak k}$ corresponding to the fundamental system.
\end{enumerate}
\end{lemma}

\begin{remark}
In Lemma~\ref{subalgebra}, the dimension of $C(\mathfrak k)$, i.e. the center of $\mathfrak k$, is 1 for case (1); 0 for cases (2) and (3), $\beta_0$ in case (3) is the highest weight of $\mathrm{ad}_{\mathfrak m^{\mathbb C}}{\mathfrak k}$ for $\theta=\theta_0$ corresponding to $\Pi'$.
\end{remark}

Let $\{\alpha_1,\alpha_2,\alpha_3,\alpha_4\}$ be a fundamental system of $\s(4)$ such that the
Dynkin diagram is
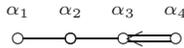
\begin{center}
\setlength{\unitlength}{0.7mm}
\begin{picture}(25,12)(0,0)
\multiput(0,0)(10,0){4}{\circle{2}}
\put(10,0){\circle{2}} \multiput(1,0)(10,0){1}{\line(1,0){8}}
\multiput(21,-0.5)(0,1){2}{\line(1,0){8}} \put(11,0){\line(1,0){8}}
\put(20.2,-1.3){$<$} \put(0,5){\makebox(0,0){$\scriptstyle \alpha_1$}}
\put(10,5){\makebox(0,0){$\scriptstyle \alpha_2$}}
\put(20,5){\makebox(0,0){$\scriptstyle \alpha_3$}}
\put(30,5){\makebox(0,0){$\scriptstyle \alpha_4$}}
\end{picture}
\end{center}
Let $\phi=2\alpha_1+2\alpha_2+2\alpha_3+\alpha_4$. By the above lemmas, the irreducible symmetric space $\Sp(4)/2\Sp(2)$ corresponds to the involution
$\theta=e^{\mathrm{ad}H}$ of $\s(4)$ defined by $$\langle H,\alpha_2\rangle=\pi\sqrt{-1};\quad \langle H,\alpha_j\rangle=0,
j=1,3,4.$$ Furthermore the Dynkin
diagram of ${2\s(2)}$ is
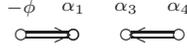
\begin{center}
\setlength{\unitlength}{0.7mm}
\begin{picture}(30,11)(0,0)
\multiput(0,0)(10,0){4}{\circle{2}}
\put(10,0){\circle{2}}
\multiput(1,-0.5)(0,1){2}{\line(1,0){8}} \put(6,-1.3){$>$}
\multiput(21,-0.5)(0,1){2}{\line(1,0){8}} \put(20.2,-1.3){$<$}
\put(0,5){\makebox(0,0){$\scriptstyle -\phi$}}
\put(10,5){\makebox(0,0){$\scriptstyle \alpha_1$}}
\put(20,5){\makebox(0,0){$\scriptstyle \alpha_3$}}
\put(30,5){\makebox(0,0){$\scriptstyle \alpha_4$}}
\end{picture}
\end{center}
Consider the involution $\tau^{1}=e^{\mathrm{ad}H_1}$ of the first $\s(2)$ with the Dynkin diagram \begin{center}
\setlength{\unitlength}{0.7mm}
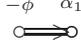
\begin{picture}(15,11)(0,0)
\multiput(0,0)(10,0){2}{\circle{2}}
\put(10,0){\circle{2}}
\multiput(1,-0.5)(0,1){2}{\line(1,0){8}} \put(6,-1.3){$>$}
\put(0,5){\makebox(0,0){$\scriptstyle -\phi$}}
\put(10,5){\makebox(0,0){$\scriptstyle \alpha_1$}}
\end{picture}
\end{center} defined by
$$\langle H_1,\alpha_1\rangle=\pi\sqrt{-1};\quad \langle
H_1,-\phi\rangle=0.$$ Let $\phi_1=-\phi+2\alpha_1=-2\alpha_2-2\alpha_3-\alpha_4$. The Dynkin diagram of $\{x\in \s(2)|\tau^{1}(x)=x\}$ is
\begin{center}
\setlength{\unitlength}{0.7mm}
\begin{picture}(15,11)(0,0)
\multiput(0,0)(10,0){2}{\circle{2}}
\put(10,0){\circle{2}}
\put(0,5){\makebox(0,0){$\scriptstyle -\phi$}}
\put(10,5){\makebox(0,0){$\scriptstyle -\phi_1$}}
\end{picture}
\end{center}
Denote by $\mathfrak h_1$ the first $\s(1)$ in the above diagram and $\mathfrak h_2$ the second $\s(1)$. Then $\mathfrak u_1=\{x\in \s(2)|\tau^{1}(x)=-x\}$ is an irreducible $\ad(\mathfrak h_1)$-module and $\s(2)=\mathfrak h_1\oplus\mathfrak h_2\oplus\mathfrak u_1$.
Similarly, consider the involution $\tau^{2}=e^{\mathrm{ad}H_2}$ of the second $\s(2)$ with the Dynkin diagram
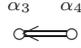
\begin{center}
\setlength{\unitlength}{0.7mm}
\begin{picture}(15,11)(0,0)
\multiput(0,0)(10,0){2}{\circle{2}}
\multiput(01,-0.5)(0,1){2}{\line(1,0){8}} \put(0.2,-1.3){$<$}
\put(0,5){\makebox(0,0){$\scriptstyle \alpha_3$}}
\put(10,5){\makebox(0,0){$\scriptstyle \alpha_4$}}
\end{picture}
\end{center}
defined by
$$\langle H_2,\alpha_3\rangle=\pi\sqrt{-1};\quad \langle H_2,\alpha_4\rangle=0.$$ Let $\phi_2=2\alpha_3+\alpha_4$. The Dynkin diagram of $\{x\in \s(2)|\tau^{2}(x)=x\}$ is
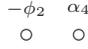
\begin{center}
\setlength{\unitlength}{0.7mm}
\begin{picture}(15,11)(0,0)
\multiput(0,0)(10,0){2}{\circle{2}}
\put(0,5){\makebox(0,0){$\scriptstyle -\phi_2$}}
\put(10,5){\makebox(0,0){$\scriptstyle \alpha_4$}}
\end{picture}
\end{center}
Denote by $\mathfrak h_3$ the first $\s(1)$ in the above diagram and $\mathfrak h_4$ the second $\s(1)$. Then $\mathfrak u_2=\{x\in \s(2)|\tau^{2}(x)=-x\}$ is an irreducible $\ad(\mathfrak h_2)$-module and $\s(2)=\mathfrak h_3\oplus\mathfrak h_4\oplus\mathfrak u_2$. Thus $\tau^{2\s(2)}=\tau^1\oplus\tau^2$ is an involution of $2\s(2)$ with the decomposition
$$2\s(2)=\mathfrak h_1\oplus\mathfrak h_2\oplus\mathfrak h_3\oplus\mathfrak h_4\oplus\mathfrak u=\mathfrak h_1\oplus\mathfrak h_2\oplus\mathfrak h_3\oplus\mathfrak h_4\oplus\mathfrak u_1\oplus \mathfrak u_2.$$
Here $\mathfrak u_1$ and $\mathfrak u_2$ are irreducible $\ad(4\s(1))$-modules. But it is unclear for $\mathfrak p$ as an $\ad(4\s(1))$-module. In order to do this, we study the extension of the involution $\tau^{2\s(2)}$ of $2\s(2)$ to $\s(4)$.

The theory on the extension of involutions of $\mathfrak k$ to $\mathfrak g$ can be found in \cite{Be2}, which is different in the method from that in \cite{Ya1}. There are also some related discussion in \cite{ChKaLi,CL,ChLi,CH}. In the following we first give the theory in general cases. Now for any involution $\tau^{\mathfrak k}$ of $\mathfrak k$, we can write
$\tau^{\mathfrak k}=\tau^{\mathfrak k}_0e^{\mathrm{ad}H^{\mathfrak k}},$ where $\tau^{\mathfrak k}_0$
is an involution on $\mathfrak k$ which keeps the Dynkin diagram of $\mathfrak k$ invariant, $H^{\mathfrak k}\in{\mathfrak t_1}$ and $\tau^{\mathfrak k}_0(H^{\mathfrak k})=H^{\mathfrak k}$.
Since $e^{\mathrm{ad}H^{\mathfrak k}}$ is an inner-automorphism, naturally we can extend $e^{\mathrm{ad}H^{\mathfrak k}}$ to an automorphism of $\mathfrak g$. Moreover,
\begin{lemma}[\cite{Ya1}]
The involution $\tau^{\mathfrak k}_0$ can be extended to an automorphism of $\mathfrak g$ if and only if $\tau^{\mathfrak k}_0$ keeps the weight system of $\mathrm{ad}_{\mathfrak m^{\mathbb C}}{\mathfrak k}$ invariant.
\end{lemma}

If $C(\mathfrak k)\not=0$, then $\dim C(\mathfrak k)=1$. Thus $\tau^{\mathfrak k}_0(Z)=Z$ or $\tau^{\mathfrak k}_0(Z)=-Z$ for any $Z\in C(\mathfrak k)$.

\begin{lemma}[\cite{Ya1}]
Assume that $C(\mathfrak k)\not=0$ and $\tau^{\mathfrak k}_0(Z)=Z$ for any $Z\in C(\mathfrak k)$. If $\tau^{\mathfrak k}$ can be extended to an automorphism of $\mathfrak g$, then $\tau^{\mathfrak k}$ can be extended to an involution of $\mathfrak g$.
\end{lemma}

\begin{lemma}[\cite{Ya1}]
Assume that $C(\mathfrak k)=0$, or $C(\mathfrak k)\not=0$ but $\tau^{\mathfrak k}_0(Z)=-Z$ for any $Z\in C(\mathfrak k)$. If $\tau$ is an automorphism of $\mathfrak g$ extending an involution $\tau^{\mathfrak k}$ of $\mathfrak k$, then $\tau^2=Id$ or $\tau^2=\theta$. Furthermore, the following conditions are equivalent:
\begin{enumerate}
   \item There exists an automorphism $\tau$ of $\mathfrak g$ extending $\tau^{\mathfrak k}$ which is an involution.
   \item Every automorphism $\tau$ of $\mathfrak g$ extending $\tau^{\mathfrak k}$ is an involution.
\end{enumerate}
\end{lemma}

\begin{lemma}[\cite{Ya1}]\label{lem2.8}
Let $\tau_0$ be the automorphism of $\mathfrak g$ extending the involution $\tau^{\mathfrak k}_0$ on $\mathfrak k$. Then $\tau_0^2=Id$ except $\mathfrak g=A_n^i$ and $n$ is even. For $e^{\mathrm{ad}H^{\mathfrak k}}$, we have:
\begin{enumerate}
 \item If $\theta_0\not=Id$, then the natural extension of $e^{\mathrm{ad}H^{\mathfrak k}}$ is an involution.
 \item Assume that $\theta_0=Id$. Let $\alpha'_{i_1},\dots,\alpha'_{i_k}$ be the roots satisfying $\langle \alpha'_j,H\rangle\not=0$. Then the natural extension of $e^{\mathrm{ad}H^{\mathfrak k}}$ is an involution if and only if $\sum_{j=i_1}^{i_k}m'_j$ is even.
\end{enumerate}
\end{lemma}

\begin{lemma}[\cite{Be2,Ya1}]\label{lem}
If $\tau$ is an involution of $\mathfrak g$ extending an involution $\tau^{\mathfrak k}$ on $\mathfrak k$, then every extension of $\tau^{\mathfrak k}$ is an involution of $\mathfrak k$, which is equivalent with $\tau$ or $\tau\theta$.
\end{lemma}

Now we go back to discuss the above case for the standard triple $(\Sp(4),2\Sp(2),4\Sp(1))$. By $(2)$ of Lemma~\ref{lem2.8}, the natural extension $\tau$ of $\tau^{2\s(2)}$ is an involution of $\s(4)$. In fact, for the fundamental system $\{\alpha_1,\alpha_2,\alpha_3,\alpha_4\}$ of $\s(4)$, the involution $\tau=e^{\mathrm{ad}H}$ of $\s(4)$ is defined by $$\langle H,\alpha_1\rangle=\langle H,\alpha_3\rangle=\pi\sqrt{-1};\quad \langle H,\alpha_2\rangle=\langle H,\alpha_4\rangle=0.$$ Clearly, $\theta\tau=\tau\theta$. That is, theorem~\ref{structure} holds.

Let ${\mathfrak p}_1=\{x\in {\mathfrak p}|\tau(x)=x\}$, and ${\mathfrak p}_2=\{x\in {\mathfrak p}|\tau(x)=-x\}$. By the above discussion, we have the following decomposition:
$$\s(4)=\mathfrak h_1\oplus\mathfrak h_2\oplus\mathfrak h_3\oplus\mathfrak h_4\oplus \mathfrak u_1\oplus \mathfrak u_2\oplus \mathfrak p_1\oplus\mathfrak p_2.$$
Consider another fundamental system $\{\alpha_2,-\alpha_1-\alpha_2,\alpha_1+\alpha_2+\alpha_3,\alpha_4\}$ of $\s(4)$ with the Dynkin diagram
\begin{center}
\setlength{\unitlength}{0.7mm}
\begin{picture}(60,12)(0,0)
\multiput(0,0)(20,0){4}{\circle{2}}
\multiput(1,0)(0,0){1}{\line(1,0){18}}
\multiput(41,-0.5)(0,1){2}{\line(1,0){18}} \put(21,0){\line(1,0){18}}
\put(40.2,-1.3){$<$} \put(0,5){\makebox(0,0){$\scriptstyle \alpha_2$}}
\put(20,5){\makebox(0,0){$\scriptstyle -\alpha_1-\alpha_2$}}
\put(40,5){\makebox(0,0){$\scriptstyle \alpha_1+\alpha_2+\alpha_3$}}
\put(60,5){\makebox(0,0){$\scriptstyle \alpha_4$}}
\end{picture}
\end{center}
The maximal root is $-\phi_1$. The Dynkin diagram of $\mathfrak h\oplus \mathfrak p_1=\{x\in \s(4)|\tau(x)=x\}$ is
\begin{center}
\setlength{\unitlength}{0.7mm}
\begin{picture}(60,11)(0,0)
\multiput(0,0)(20,0){4}{\circle{2}}
\multiput(1,-0.5)(0,1){2}{\line(1,0){18}} \put(16,-1.3){$>$}
\multiput(41,-0.5)(0,1){2}{\line(1,0){18}} \put(40.2,-1.3){$<$}
\put(0,5){\makebox(0,0){$\scriptstyle -\phi_1$}}
\put(20,5){\makebox(0,0){$\scriptstyle \alpha_2$}}
\put(40,5){\makebox(0,0){$\scriptstyle \alpha_1+\alpha_2+\alpha_3$}}
\put(60,5){\makebox(0,0){$\scriptstyle \alpha_4$}}
\end{picture}
\end{center}
Consider the restriction of the involution $\theta$ of $\s(4)$ on $\mathfrak h\oplus \mathfrak p_1$. Similar to the discussion for $\mathfrak u$, $$\mathfrak p_1=\mathfrak p_1^1\oplus \mathfrak p_1^2$$ as the direct sum of irreducible $\ad(4\s(1))$-modules. For $\mathfrak p_2$, consider the fundamental system $\{\alpha_1+\alpha_2,-\alpha_2,\alpha_2+\alpha_3,\alpha_4\}$ of $\s(4)$ with the Dynkin diagram
\begin{center}
\setlength{\unitlength}{0.7mm}
\begin{picture}(60,12)(0,0)
\multiput(0,0)(20,0){4}{\circle{2}}
\multiput(1,0)(0,0){1}{\line(1,0){18}}
\multiput(41,-0.5)(0,1){2}{\line(1,0){18}} \put(21,0){\line(1,0){18}}
\put(40.2,-1.3){$<$} \put(0,5){\makebox(0,0){$\scriptstyle \alpha_1+\alpha_2$}}
\put(20,5){\makebox(0,0){$\scriptstyle -\alpha_2$}}
\put(40,5){\makebox(0,0){$\scriptstyle \alpha_2+\alpha_3$}}
\put(60,5){\makebox(0,0){$\scriptstyle \alpha_4$}}
\end{picture}
\end{center}
The maximal root is $\phi$. The Dynkin diagram of $\mathfrak h\oplus \mathfrak p_2=\{x\in \s(4)|\theta\tau(x)=x\}$ is
\begin{center}
\setlength{\unitlength}{0.7mm}
\begin{picture}(60,11)(0,0)
\multiput(0,0)(20,0){4}{\circle{2}}
\multiput(1,-0.5)(0,1){2}{\line(1,0){18}} \put(16,-1.3){$>$}
\multiput(41,-0.5)(0,1){2}{\line(1,0){18}} \put(40.2,-1.3){$<$}
\put(0,5){\makebox(0,0){$\scriptstyle -\phi$}}
\put(20,5){\makebox(0,0){$\scriptstyle \alpha_1+\alpha_2$}}
\put(40,5){\makebox(0,0){$\scriptstyle \alpha_2+\alpha_3$}}
\put(60,5){\makebox(0,0){$\scriptstyle \alpha_4$}}
\end{picture}
\end{center}
Consider the restriction of the involution $\theta$ of $\s(4)$ on $\mathfrak h\oplus \mathfrak p_2$. Similar to the discussion for $\mathfrak u$, $$\mathfrak p_2=\mathfrak p_2^1\oplus \mathfrak p_2^2$$ as the direct sum of irreducible $\ad(4\s(1))$-modules. Up to now, we have the decomposition
$$\s(4)=\mathfrak h_1\oplus\mathfrak h_2\oplus\mathfrak h_3\oplus\mathfrak h_4\oplus \mathfrak u_1\oplus \mathfrak u_2\oplus \mathfrak p_1^1\oplus \mathfrak p_1^2\oplus \mathfrak p_2^1\oplus \mathfrak p_2^2$$
as a direct sum of irreducible $\ad(4\s(1))$-modules.

\section{Non-naturally reductive Einstein metrics on $\Sp(4)$ which are $\Ad(4\Sp(1))$-invariant}\label{general}
For the standard triple $(\Sp(4),2\Sp(2),4\Sp(1)$, we have the following $B$-orthogonal decomposition of $\s(4)$ as irreducible $\ad(4\s(1))$-modules: $$\s(4)=\mathfrak h_1\oplus\mathfrak h_2\oplus\mathfrak h_3\oplus\mathfrak h_4\oplus \mathfrak u_1\oplus \mathfrak u_2\oplus \mathfrak p_1^1\oplus \mathfrak p_1^2\oplus \mathfrak p_2^1\oplus \mathfrak p_2^2.$$
Consider the left-invariant metric $\langle\cdot,\cdot\rangle$ on $\Sp(4)$ which is $\Ad(H)$-invariant of the form
\begin{equation}\label{met-fm}
x_1B|_{\mathfrak h_1}\oplus x_2B|_{\mathfrak h_2}\oplus x_3B|_{\mathfrak h_3}\oplus x_4B|_{\mathfrak h_4}\oplus x_5B|_{\mathfrak u_1}\oplus x_6B|_{\mathfrak u_2}\oplus x_7B|_{\mathfrak p_1^1}\oplus x_8B|_{\mathfrak p_1^2}\oplus x_9B|_{\mathfrak p_2^1}\oplus x_{10}B|_{\mathfrak p_2^2},\end{equation}
In \cite{YD1}, the authors show that there are no non-naturally reductive Einstein metrics which satisfy $$x_1=x_2=x_3=x_4, \quad x_5=x_6, \quad x_7=x_8=x_9=x_{10}.$$
In the following, we consider the metrics on $\Sp(4)$ of the form~\eqref{met-fm} satisfying $$x_3=x_4=x_6, \quad x_7=x_{10}, \quad x_8=x_9.$$
That is, we can consider the following $B$-orthogonal decomposition
$$\s(4)=\mathfrak h_1\oplus\mathfrak h_2\oplus\{\mathfrak h_3\oplus\mathfrak h_4\oplus \mathfrak u_2\}\oplus \mathfrak u_1\oplus \{\mathfrak p_1^1\oplus \mathfrak p_2^2\}\oplus \{\mathfrak p_1^2\oplus \mathfrak p_2^1\}.$$
The above decomposition corresponds to another involution pair $(\theta,\sigma)$ satisfying $\theta\sigma=\sigma\theta$, where $\theta$ is the involution in section~\ref{sec2} and $\sigma=e^{\mathrm{ad}H}$ is the involution defined by $$\langle H,\alpha_1\rangle=\pi\sqrt{-1};\quad \langle H,\alpha_j\rangle=0, j=2,3,4.$$
Here $\sigma$ is the natural extension of the involution $e^{\mathrm{ad}H'}$ on $\mathfrak k$ defined by
$$\langle H',\alpha_1\rangle=\pi\sqrt{-1};\quad \langle H',-\phi\rangle=\langle H',\alpha_3\rangle=\langle H',\alpha_4\rangle=0.$$
We can check that $\mathfrak h_1\oplus\mathfrak h_2\oplus\{\mathfrak h_3\oplus\mathfrak h_4\oplus \mathfrak u_2\}=\s(1)\oplus\s(1)\oplus\s(2)=\{x\in\s(4)|\theta(x)=\sigma(x)=x\}$, and $\Sp(4)/\Sp(1)\times\Sp(1)\times \Sp(2)$ is a generalized Wallach space.

In general, for any $1\leq k_1\leq k_2\leq k_3$, $\Sp(k_1+k_2+k_3)/\Sp(k_1)\times \Sp(k_2)\times \Sp(k_3)$ is a generalized Wallach space \cite{ChKaLi,Ni16}. Let $\{\alpha_1,\alpha_2,\cdots,\alpha_{k_1+k_2+k_3}\}$ be a fundamental system of $\s(k_1+k_2+k_3)$ with the Dynkin diagram
\begin{center}
\setlength{\unitlength}{0.7mm}
\begin{picture}(25,12)(0,0)
\multiput(0,0)(10,0){4}{\circle{2}}
\put(10,0){\circle{2}} \multiput(1,0)(10,0){1}{\line(1,0){8}}
\multiput(21,-0.5)(0,1){2}{\line(1,0){8}}
\put(11,0){\line(1,0){1}}
\put(13,0){\line(1,0){1}}
\put(15,0){\line(1,0){1}}
\put(17,0){\line(1,0){1}}
\put(20.2,-1.3){$<$} \put(0,5){\makebox(0,0){$\scriptstyle \alpha_1$}}
\put(10,5){\makebox(0,0){$\scriptstyle \alpha_2$}}
\put(30,5){\makebox(0,0){$\scriptstyle \alpha_{k_1+k_2+k_3}$}}
\end{picture}
\end{center}
The generalized Wallach space $\Sp(k_1+k_2+k_3)/\Sp(k_1)\times \Sp(k_2)\times \Sp(k_3)$ corresponds to the involution pair $(\theta,\tau)$ of $\s(k_1+k_2+k_3)$, where $\theta=e^{\mathrm{ad}H}$ determined by $$\langle H,\alpha_{k_1+k_2}\rangle=\pi\sqrt{-1};\quad \langle H,\alpha_j\rangle=0,
j\not=k_1+k_2,$$
and $\tau=e^{\mathrm{ad}H'}$ determined by $$\langle H',\alpha_{k_1}\rangle=\pi\sqrt{-1};\quad \langle H',\alpha_j\rangle=0,
j\not=k_1.$$
Then we have the decomposition of $\s(k_1+k_2+k_3)$ as irreducible $\ad(\s(k_1)\oplus \s(k_2)\oplus \s(k_3))$-modules:
$$\s(k_1+k_2+k_3)=\s(k_1)\oplus \s(k_2)\oplus \s(k_3)\oplus \mathfrak m_1\oplus \mathfrak m_2\oplus \mathfrak m_3,$$
where $\mathfrak m_1=\{x|\theta(x)=-\tau(x)=x\}$, $\mathfrak m_2=\{x|\theta(x)=\tau(x)=-x\}$, $\mathfrak m_3=\{x|\theta(x)=-\tau(x)=-x\}$. Consider the left-invariant metrics on $\Sp(k_1+k_2+k_3)$ determined by the $\Ad(\Sp(k_1)\times \Sp(k_2)\times \Sp(k_3))$-invariant inner product on $\s(k_1+k_2+k_3)$ given by
\begin{equation}\label{metricsp}
\langle\cdot,\cdot\rangle=y_1B|_{\s(k_1)}\oplus y_2B|_{\s(k_2)}\oplus y_3B|_{\s(k_3)}\oplus y_4B|_{\mathfrak m_1}\oplus y_5B|_{\mathfrak m_2}\oplus y_6B|_{\mathfrak m_3}.
\end{equation}
By the theory given in Lemma~\ref{nr},
\begin{lemma}[\cite{ArSaSt}]\label{naturallyred}
If a left-invariant metric of the form~\eqref{metricsp} is naturally reductive with respect to $\Sp(k_1+k_2+k_3)\times L$ for some closed subgroup $L$ of $\Sp(k_1+k_2+k_3)$, then one of the following holds:
\begin{enumerate}
 \item $y_1=y_2=y_4$, $y_5=y_6$.
 \item $y_2=y_3=y_6$, $y_4=y_5$.
 \item $y_1=y_3=y_5$, $y_4=y_6$.
 \item $y_4=y_5=y_6$.
\end{enumerate}
Conversely, if one of the above conditions holds, then the metric of the form~~\eqref{metricsp} is naturally reductive with respect to $\Sp(k_1+k_2+k_3)\times L$ for some closed subgroup $L$ of $\Sp(k_1+k_2+k_3)$.
\end{lemma}
Using the formula in Lemma~\ref{formula1}, we have the following lemma.
\begin{lemma}[\cite{ArSaSt}]
The components of the Ricci tensor for the metric of the form~\eqref{metricsp} are: \begin{eqnarray*}
   && r_{\Sp(1)}=\frac{k_1+1}{4(n+1)y_1}+\frac{k_2y_1}{4(n+1)y_4^2}+\frac{k_3y_1}{4(n+1)y_5^2}, \\
    && r_{\Sp(2)}=\frac{k_2+1}{4(n+1)y_2}+\frac{k_1y_2}{4(n+1)y_4^2}+\frac{k_3y_2}{4(n+1)y_6^2},\\
 && r_{\Sp(3)}=\frac{k_3+1}{4(n+1)y_3}+\frac{k_1y_3}{4(n+1)y_5^2}+\frac{k_2y_3}{4(n+1)y_6^2}, \\
 &&
 r_{\mathfrak m_1}=\frac{1}{2y_4}+\frac{k_3}{4(n+1)}\left(\frac{y_4}{y_5y_6}-\frac{y_5}{y_4y_6}-\frac{y_6}{y_4y_5}\right)
 -\frac{(2k_1+1)y_1}{8(n+1)y_4^2}-\frac{(2k_2+1)y_2}{8(n+1)y_4^2}, \\
 &&
 r_{\mathfrak m_2}=\frac{1}{2y_5}+\frac{k_2}{4(n+1)}\left(\frac{y_5}{y_4y_6}-\frac{y_4}{y_5y_6}-\frac{y_6}{y_4y_5}\right)
 -\frac{(2k_1+1)y_1}{8(n+1)y_5^2}-\frac{(2k_3+1)y_3}{8(n+1)y_5^2}, \\&&
 r_{\mathfrak m_3}=\frac{1}{2y_6}+\frac{k_1}{4(n+1)}\left(\frac{y_6}{y_4y_5}-\frac{y_4}{y_5y_6}-\frac{y_5}{y_4y_6}\right)
 -\frac{(2k_2+1)y_2}{8(n+1)y_6^2}-\frac{(2k_3+1)y_3}{8(n+1)y_6^2}.
\end{eqnarray*}
Here $n=k_1+k_2+k_3$.
\end{lemma}
Furthermore, for $k_1=k_2=1$, $k_3=2$, we have the following lemma.
\begin{lemma}[\cite{ArSaSt}]\label{lemma4.3}
The Lie group $\Sp(4)$ admits three $\Ad(\Sp(1)\times\Sp(1)\times\Sp(2))$-invariant Einstein metrics of the form~\eqref{metricsp} which are non-naturally reductive up to isometry and scalar. The solutions for the non-naturally reductive metrics are given by
\begin{enumerate}
 \item $(y_1, y_2, y_3, y_4, y_5, y_6)\approx
(0.114935, 0.114935, 0.180564, 0.508812, 0.326608, 0.326608)$,
 \item $(y_1, y_2, y_3, y_4, y_5, y_6)\approx
(0.116403, 0.116403, 0.169957, 0.310184, 0.380445, 0.380445)$,
 \item $(y_1, y_2, y_3, y_4, y_5, y_6)\approx
(0.117632, 0.131837, 0.170185, 0.241674, 0.489011, 0.320149)$.
\end{enumerate}
\end{lemma}

Now we are at the point to prove Theorem~\ref{admitmetric}. In fact, for the values of $y_i$ in Lemma~\ref{lemma4.3}, every metric of the form~\eqref{met-fm} satisfying $$x_1=y_1,x_2=y_2,x_3=x_4=x_6=y_3,x_5=y_4,x_7=x_{10}=y_5,x_8=x_9=y_6$$
determines a non-naturally reductive Einstein metric on $\Sp(4)$ which is $\Ad(4\Sp(1))$-invariant.

\section{The decomposition of $\s(2n_1n_2)$ as a direct sum of irreducible $\ad(2n_2\s(n_1))$-modules.}
In the following, we first prove Theorem~\ref{structure1}. In fact, let $\{\alpha_1,\alpha_2,\cdots,\alpha_{2n_1n_2}\}$ be a fundamental system of $\s(2n_1n_2)$ with the Dynkin diagram
\begin{center}
\setlength{\unitlength}{0.7mm}
\begin{picture}(25,12)(0,0)
\multiput(0,0)(10,0){4}{\circle{2}}
\put(10,0){\circle{2}} \multiput(1,0)(10,0){1}{\line(1,0){8}}
\multiput(21,-0.5)(0,1){2}{\line(1,0){8}}
\put(11,0){\line(1,0){1}}
\put(13,0){\line(1,0){1}}
\put(15,0){\line(1,0){1}}
\put(17,0){\line(1,0){1}}
\put(20.2,-1.3){$<$} \put(0,5){\makebox(0,0){$\scriptstyle \alpha_1$}}
\put(10,5){\makebox(0,0){$\scriptstyle \alpha_2$}}
\put(30,5){\makebox(0,0){$\scriptstyle \alpha_{2n_1n_2}$}}
\end{picture}
\end{center}
For any $1\leq i\leq 2n_2-1$, let $\theta_i$ be the involution of $\s(2n_1n_2)$ defined by $\theta_i=e^{\mathrm{ad}H_i}$ where $H_i$ satisfies $$\langle H_i,\alpha_{n_1\,i}\rangle=\pi\sqrt{-1};\quad \langle H_i,\alpha_j\rangle=0,
j\not=n_1\,i.$$
Clearly, $\theta_i\theta_j=\theta_j\theta_i$. Then we have the following decomposition of $\s(2n_1n_2)$:
$$\s(2n_1n_2)=\sum_{j_i=1 \text{ or } -1 }\mathfrak m{(j_1,\cdots,j_{2n_2-1})},$$
Here $\mathfrak m{(j_1,\cdots,j_{2n_2-1})}=\{x\in \s(2n_1n_2)|\theta_i(x)=j_ix\}$, and in particular $$\mathfrak m{(1,\cdots,1)}=\s(n_1)^1\oplus\cdots\oplus\s(n_1)^{2n_2}.$$ Clearly $j_i=1$ or $-1$. Define $\Phi=\{(j_1,\cdots,j_{2n_2-1})|\mathfrak m{(j_1,\cdots,j_{2n_2-1})}\not=0, (j_1,\cdots,j_{2n_2-1})\not=(1,1,\cdots,1)\}$. Then we have
$$\s(2n_1n_2)=\s(n_1)^1\oplus\cdots\oplus\s(n_1)^{2n_2}\oplus\sum_{(j_1,\cdots,j_{2n_2-1})\in\Phi}\mathfrak m{(j_1,\cdots,j_{2n_2-1})},$$ which is a direct sum of irreducible $\ad(2n_2\s(n_1))$-modules.

\begin{remark}\label{rem}
The standard triple $(\Sp(4),2\Sp(2),4\Sp(1))$ corresponds to the above case with $n_1=1$, $n_2=2$. Then $$\s(4)=\s(n_1)^1\oplus\cdots\oplus\s(n_1)^{4}\oplus\sum_{(j_1,j_2,j_3)\in\Phi}\mathfrak m{(j_1,j_2,j_3)}.$$
Here there are possibly 11 irreducible $\ad(4\s(1))$-modules: $\s(n_1)^1$, $\s(n_1)^2$, $\s(n_1)^{3}$, $\s(n_1)^{4}$, $\mathfrak m(1,1,-1)$, $\mathfrak m(1,-1,1)$, $\mathfrak m(-1,1,1)$, $\mathfrak m(1,-1,-1)$, $\mathfrak m(-1,1,-1)$, $\mathfrak m(-1,-1,1)$, and $\mathfrak m(-1,-1,-1)$. But we prove in section~\ref{sec2} that $\s(4)$ is the direct sum of 10 irreducible $\ad(4\s(1))$-modules. In fact, it is easy to check that $\mathfrak m{(-1,1,-1)}=0$.
\end{remark}

Consider the left-invariant metrics on $\Sp(2n_1n_2)$ determined by the $\Ad(2n_2\Sp(n_1))$-invariant inner product on $\s(2n_1n_2)$ given by
\begin{equation}\label{metricsp1}
\langle\cdot,\cdot\rangle=\sum_{i=1}^{2n_2}y_iB|_{\s(n_1)^i} \oplus \sum_{(j_1,\cdots,j_{2n_2-1})\in\Phi} y{(j_1,\cdots,j_{2n_2-1})}B|_{\mathfrak m_(j_1,\cdots,j_{2n_2-1})}.
\end{equation}
In \cite{YD1}, the authors show that there are non-naturally reductive Einstein metrics on $\Sp(2n_1n_2)$ for $n_1n_2>2$ which satisfy
\begin{enumerate}
  \item $y_1=y_2=\cdots=y_{2n_2}$, and
  \item $y{(j_1,\cdots,j_{2n_2-1})}=m$ if $(j_1,\cdots,j_{n_2}=-1,\cdots,j_{2n_2-1})\in\Phi$, and
  \item $y{(j_1,\cdots,j_{2n_2-1})}=n$ if $(j_1,\cdots,j_{n_2}=1,\cdots,j_{2n_2-1})\in\Phi$.
\end{enumerate}
We will consider the metric of the form~\eqref{metricsp1} satisfying
\begin{enumerate}
  \item $y_3=\cdots=y_{2n_2}=y{(j_1,\cdots,j_{2n_2-1})}$ for any $(j_1=1,j_2=1,j_3,\cdots,j_{2n_2-1})\in\Phi$, and
  \item all $y{(j_1,\cdots,j_{2n_2-1})}$ satisfying $(j_1=1,j_2=-1,j_3,\cdots,j_{2n_2-1})\in\Phi$ are equivalent, and
  \item all $y{(j_1,\cdots,j_{2n_2-1})}$ satisfying $(j_1=-1,j_2=1,j_3,\cdots,j_{2n_2-1})\in\Phi$ are equivalent, and
  \item all $y{(j_1,\cdots,j_{2n_2-1})}$ satisfying $(j_1=-1,j_2=-1,j_3,\cdots,j_{2n_2-1})\in\Phi$ are equivalent.
\end{enumerate}
The above metric corresponds to the decomposition of $\s(2n_1n_2)$ under the involution pair $(\theta_1,\theta_2)$ which are $\Ad(\Sp(n_1)\times\Sp(n_1)\times\Sp(2(n_2-1)n_1))$-invariant. Generally, we will study non-naturally reductive Einstein metrics on $\Sp(2k+l)$ which are $\Ad(\Sp(k)\times\Sp(k)\times\Sp(l))$-invariant.

\section{Non-naturally reductive Einstein metrics on $\Sp(2k+l)$}
In this section, we give non-naturally reductive Einstein metrics on $\Sp(2k+l)$ which are $\Ad(\Sp(k)\times\Sp(k)\times\Sp(l))$-invariant.
That is, we take $k=k_1=k_2$ and $l=k_3$ in the formulae given in Section~\ref{general}, and study the metric with $y_1=y_2$ and $y_5=y_6=1$. The homogeneous Einstein equations are equivalent to the following system of equations:
\[\left\{\begin{aligned}
f_1&\triangleq-2\,k{y_{{3}}}^{2}y_{{2}}{y_{{4}}}^{2}+l{y_{{2}}}^{2}y_{{3}}{y_{{4}}}^
{2}+k{y_{{2}}}^{2}y_{{3}}+ky_{{3}}{y_{{4}}}^{2}-ly_{{2}}{y_{{4}}}^{2}-
y_{{2}}{y_{{4}}}^{2}+y_{{3}}{y_{{4}}}^{2}=0,\\
f_2&\triangleq l{y_{{2}}}^{2}{y_{{4}}}^{2}-ly_{{2}}{y_{{4}}}^{3}+3\,k{y_{{2}}}^{2}-4
\,y_{{2}}y_{{4}}k+k{y_{{4}}}^{2}+{y_{{2}}}^{2}-2\,y_{{2}}y_{{4}}+{y_{{
4}}}^{2}=0,\\
f_3&\triangleq 2\,y_{{2}}y_{{3}}k+4\,k{y_{{3}}}^{2}+2\,y_{{3}}y_{{4}}k+2\,{y_{{3}}}^{
2}l-8\,ky_{{3}}-4\,y_{{3}}l+y_{{2}}y_{{3}}+{y_{{3}}}^{2}+2\,l-4\,y_{{3
}}+2=0.
\end{aligned}\right.\]
In particular,
$$f_2= \left( y_{{2}}-y_{{4}} \right)  \left( ly_{{2}}{y_{{4}}}^{2}+3\,ky_{{2}}-ky_{{4}}+y_{{2}}-y_{{4}} \right).$$
If $y_2=y_4$, by Lemma~\ref{naturallyred}, the metrics are naturally reductive. In order to get non-naturally reductive Einstein metrics, we assume that $y_2\neq y_4$. Thus $ly_{{2}}{y_{{4}}}^{2}+3\,ky_{{
2}}-ky_{{4}}+y_{{2}}-y_{{4}}=0$. That is, $$y_2={\frac { \left( k+1 \right) y_{{4}}}{l{y_{{4}}}^{2}+3\,k+1}}.$$ Substituting it into $f_1$ and $f_3$, we have
\[\left\{\begin{aligned}
g_1&\triangleq 2\,{k}^{2}l{y_{{3}}}^{2}{y_{{4}}}^{3}-k{l}^{2}y_{{3}}{y_{{4}}}^{4}+2\,
kl{y_{{3}}}^{2}{y_{{4}}}^{3}-{l}^{2}y_{{3}}{y_{{4}}}^{4}+6\,{k}^{3}{y_
{{3}}}^{2}y_{{4}}-7\,{k}^{2}ly_{{3}}{y_{{4}}}^{2}+k{l}^{2}{y_{{4}}}^{3
}+8\,{k}^{2}{y_{{3}}}^{2}y_{{4}}\\&-10\,kly_{{3}}{y_{{4}}}^{2}+kl{y_{{4}}
}^{3}+{l}^{2}{y_{{4}}}^{3}-10\,{k}^{3}y_{{3}}+3\,{k}^{2}ly_{{4}}+2\,k{
y_{{3}}}^{2}y_{{4}}-3\,ly_{{3}}{y_{{4}}}^{2}+l{y_{{4}}}^{3}-17\,{k}^{2
}y_{{3}}+3\,{k}^{2}y_{{4}}\\&+4\,kly_{{4}}-8\,ky_{{3}}+4\,ky_{{4}}+ly_{{4
}}-y_{{3}}+y_{{4}}=0,\\
g_3&\triangleq 4\,kl{y_{{3}}}^{2}{y_{{4}}}^{2}+2\,kly_{{3}}{y_{{4}}}^{3}+2\,{l}^{2}{y
_{{3}}}^{2}{y_{{4}}}^{2}-8\,kly_{{3}}{y_{{4}}}^{2}-4\,{l}^{2}y_{{3}}{y
_{{4}}}^{2}+l{y_{{3}}}^{2}{y_{{4}}}^{2}+12\,{k}^{2}{y_{{3}}}^{2}+8\,y_
{{3}}y_{{4}}{k}^{2}\\&+6\,kl{y_{{3}}}^{2}+2\,{l}^{2}{y_{{4}}}^{2}-4\,ly_{
{3}}{y_{{4}}}^{2}-24\,{k}^{2}y_{{3}}-12\,kly_{{3}}+7\,k{y_{{3}}}^{2}+5
\,y_{{3}}y_{{4}}k+2\,{y_{{3}}}^{2}l+2\,l{y_{{4}}}^{2}+6\,kl\\&-20\,ky_{{3
}}-4\,y_{{3}}l+{y_{{3}}}^{2}+y_{{3}}y_{{4}}+6\,k+2\,l-4\,y_{{3}}+2=0.
\end{aligned}\right.\]
Consider the polynomial ring $R=\mathbb{Q}[z, y_3, y_4]$ and the ideal $I$ generated by $\{g_1, g_3, z y_3 y_4-1\}$, and take a lexicographic order $>$ with $z>y_3>y_4$ for a monomial ordering on $R$. By the help of computer, we get two polynomials in the Gr\"{o}bner basis of the ideal $I$. One is \\
$h(y_4)=2\,{l}^{2}( l+k ) ( 4\,{k}^{2}+4\,kl+2\,{l}^{2}+l
  )y_4^8 \\
  \indent -4\,{l}^{2}  ( 2\,k+l+1  )   ( 8\,{k}^{2}+8\,kl+2\,{l}^{2
}+l  ) y_4^7\\
\indent +l  ( 64\,{k}^{4}+304\,{k}^{3}l+284\,{k}^{2}{l}^{2}+88\,k{l}^{3}+4
\,{l}^{4}+40\,{k}^{3}+238\,{k}^{2}l+162\,k{l}^{2}+32\,{l}^{3}+8\,{k}^{
2}+51\,kl+19\,{l}^{2}+2\,l  ) y_4^6\\
\indent -8\,l  ( 2\,k+l+1  )   ( 28\,{k}^{3}+40\,{k}^{2}l+10\,k{l
}^{2}+14\,{k}^{2}+21\,kl+4\,{l}^{2}+2\,k+2\,l  ) y_4^5\\
\indent +(128\,{k}^{5}+1376\,{k}^{4}l+1472\,{k}^{3}{l}^{2}+468\,{k}^{2}{l}^{3}+
24\,k{l}^{4}+160\,{k}^{4}+1684\,{k}^{3}l+1372\,{k}^{2}{l}^{2}+336\,k{l
}^{3}+8\,{l}^{4}\\
\indent +82\,{k}^{3}+713\,{k}^{2}l+406\,k{l}^{2}+60\,{l}^{3}+
20\,{k}^{2}+118\,kl+38\,{l}^{2}+2\,k+5\,l )y_4^4\\
\indent -4\,  ( 2\,k+l+1  )   ( 96\,{k}^{4}+248\,{k}^{3}l+62\,{k}
^{2}{l}^{2}+92\,{k}^{3}+215\,{k}^{2}l+46\,k{l}^{2}+32\,{k}^{2}+55\,kl+
8\,{l}^{2}+4\,k+4\,l  )y_4^3\\
\indent +(1792\,{k}^{5}+2432\,{k}^{4}l+848\,{k}^{3}{l}^{2}+36\,{k}^{2}{l}^{3}+
2848\,{k}^{4}+3100\,{k}^{3}l+896\,{k}^{2}{l}^{2}+24\,k{l}^{3}+1714\,{k
}^{3}+1427\,{k}^{2}l \\
\indent +304\,k{l}^{2}+4\,{l}^{3}+478\,{k}^{2}+278\,kl+32
\,{l}^{2}+58\,k+19\,l+2)y_4^2\\
\indent -4\,  ( 5\,k+1  ) ( 2\,k+1  )   ( 3\,k+1)( 2\,k+l+1  )( 8\,k+2\,l+1  ) y_4\\
\indent
+2\,  ( 2\,k+1  ) ^{2}  ( 5\,k+1  ) ^{2}  ( 4\,k+
2\,l+1  )$,\\
the other is \\
$h(y_3,y_4)=-4\,{l}^{2}  ( l+k  )   ( 4\,{k}^{2}+4\,kl+2\,{l}^{2}+l
  )   ( {k}^{2}+4\,kl+2\,{l}^{2}+k+l  ) y_4^7\\
  \indent +8\,{l}^{2}  ( 2\,k+l+1  )   ( {k}^{2}+4\,kl+2\,{l}^{2}+k+
l  )   ( 8\,{k}^{2}+8\,kl+2\,{l}^{2}+l  ) y_4^6\\
\indent -2\,l  ( 64\,{k}^{6}+480\,{k}^{5}l+1428\,{k}^{4}{l}^{2}+1632\,{k}^
{3}{l}^{3}+824\,{k}^{2}{l}^{4}+172\,k{l}^{5}+8\,{l}^{6}+104\,{k}^{5}+
730\,{k}^{4}l+1682\,{k}^{3}{l}^{2}\\
\indent +1428\,{k}^{2}{l}^{3}+504\,k{l}^{4}+
64\,{l}^{5}+48\,{k}^{4}+357\,{k}^{3}l+622\,{k}^{2}{l}^{2}+355\,k{l}^{3
}+66\,{l}^{4}+8\,{k}^{3}+61\,{k}^{2}l+77\,k{l}^{2}\\
\indent +22\,{l}^{3}+2\,kl+2
\,{l}^{2}  ) y_4^5\\
\indent +4\,l  ( 2\,k+l+1  )   ( 112\,{k}^{5}+448\,{k}^{4}l+664\,{
k}^{3}{l}^{2}+360\,{k}^{2}{l}^{3}+60\,k{l}^{4}+168\,{k}^{4}+508\,{k}^{
3}l+556\,{k}^{2}{l}^{2}+228\,k{l}^{3}\\
\indent +28\,{l}^{4}+64\,{k}^{3}+172\,{k}
^{2}l+131\,k{l}^{2}+28\,{l}^{3}+8\,{k}^{2}+16\,kl+7\,{l}^{2}  ) y_4^4\\
\indent -(256\,{k}^{7}-2976\,{k}^{6}l-8944\,{k}^{5}{l}^{2}-10696\,{k}^{4}{l}^{3
}-5680\,{k}^{3}{l}^{4}-1224\,{k}^{2}{l}^{5}-56\,k{l}^{6}-576\,{k}^{6}-
6656\,{k}^{5}l-15528\,{k}^{4}{l}^{2}\\
\indent -14480\,{k}^{3}{l}^{3}-5996\,{k}^{
2}{l}^{4}-988\,k{l}^{5}-24\,{l}^{6}-484\,{k}^{5}-5282\,{k}^{4}l-9632\,
{k}^{3}{l}^{2}-6830\,{k}^{2}{l}^{3}-2006\,k{l}^{4}-196\,{l}^{5}\\
\indent -204\,{
k}^{4}-1882\,{k}^{3}l-2625\,{k}^{2}{l}^{2}-1293\,k{l}^{3}-210\,{l}^{4}
-44\,{k}^{3}-294\,{k}^{2}l-293\,k{l}^{2}-75\,{l}^{3}-4\,{k}^{2}-14\,kl
-8\,{l}^{2})y_4^3\\
\indent +4\,  ( 2\,k+l+1  )   ( 192\,{k}^{6}+624\,{k}^{5}l+1052\,{
k}^{4}{l}^{2}+648\,{k}^{3}{l}^{3}+108\,{k}^{2}{l}^{4}+376\,{k}^{5}+
1166\,{k}^{4}l+1412\,{k}^{3}{l}^{2}+684\,{k}^{2}{l}^{3}\\ \indent +96\,k{l}^{4}+
248\,{k}^{4}+688\,{k}^{3}l+631\,{k}^{2}{l}^{2}+216\,k{l}^{3}+20\,{l}^{
4}+72\,{k}^{3}+158\,{k}^{2}l+108\,k{l}^{2}+20\,{l}^{3}+8\,{k}^{2}+12\,
kl+5\,{l}^{2}  ) y_4^2\\
\indent -(2944\,{k}^{7}-5920\,{k}^{6}l-6800\,{k}^{5}{l}^{2}-4104\,{k}^{4}{l}^{3
}-960\,{k}^{3}{l}^{4}-24\,{k}^{2}{l}^{5}-7632\,{k}^{6}-13816\,{k}^{5}l
-13320\,{k}^{4}{l}^{2} \\
\indent -6568\,{k}^{3}{l}^{3}-1284\,{k}^{2}{l}^{4}-32\,k
{l}^{5}-7480\,{k}^{5}-11750\,{k}^{4}l-9368\,{k}^{3}{l}^{2}-3574\,{k}^{
2}{l}^{3}-512\,k{l}^{4}-8\,{l}^{5}-3546\,{k}^{4}\\
\indent -4678\,{k}^{3}l-2943\,
{k}^{2}{l}^{2}-784\,k{l}^{3}-60\,{l}^{4}-838\,{k}^{3}-893\,{k}^{2}l-
400\,k{l}^{2}-58\,{l}^{3}-86\,{k}^{2}-70\,kl-17\,{l}^{2}-2\,k-l)y_4\\
\indent +  ( 2\,l+1  )   ( 5\,k+1  )   ( 2\,k+1  )
  ( k+1  )   ( 4\,k+2\,l+1  )   ( 6\,{k}^{2}+4\,
kl+2\,{l}^{2}+2\,k+l  ) y_3\\
\indent +16\,k  ( 5\,k+1  )   ( 2\,k+1  )   ( k+1
  )   ( 2\,k+l+1  )   ( 6\,{k}^{2}+4\,kl+2\,{l}^{2}+
2\,k+l  ) $.\\
It is easy to see that there exists a real number $s$ such that $h(s,t)=0$, where $t$ is a real number such that $h(t)=0$. Furthermore, we can prove that $s>0$. In fact, we take the lexicographic order $>$ with $z>y_4>y_3$ for a monomial ordering on $R$. By the help of computer, the following polynomial $t(y_3)$ is contained in the Gr\"{o}bner basis for the ideal $I$:\\
$t(y_3)= ( 4\,k+2\,l+1 )  ( 4\,{k}^{2}+4\,kl+2\,{l}^{2}+l
 )  ( 6\,{k}^{2}+4\,kl+2\,{l}^{2}+2\,k+l ) ^{2}y_3^8\\
 \indent -16\, ( 2\,k+l+1 )  ( 6\,{k}^{2}+4\,kl+2\,{l}^{2}+2\,k
+l )  ( 6\,{k}^{4}+24\,{k}^{3}l+28\,{k}^{2}{l}^{2}+16\,k{l}
^{3}+4\,{l}^{4}+2\,{k}^{3}+10\,{k}^{2}l \\
\indent +10\,k{l}^{2}+4\,{l}^{3}+kl+{l}
^{2} ) y_3^7\\
\indent +(640\,{k}^{7}+7712\,{k}^{6}l+21968\,{k}^{5}{l}^{2}+30576\,{k}^{4}{l}^{3
}+24804\,{k}^{3}{l}^{4}+12288\,{k}^{2}{l}^{5}+3520\,k{l}^{6}+448\,{l}^
{7}+2000\,{k}^{6}\\
\indent +13528\,{k}^{5}l+31660\,{k}^{4}{l}^{2}+36532\,{k}^{3}
{l}^{3}+23504\,{k}^{2}{l}^{4}+8416\,k{l}^{5}+1328\,{l}^{6}+1408\,{k}^{
5}+7818\,{k}^{4}l+15881\,{k}^{3}{l}^{2} \\
\indent +14848\,{k}^{2}{l}^{3}+7012\,k{
l}^{4}+1432\,{l}^{5}+368\,{k}^{4}+1924\,{k}^{3}l+3456\,{k}^{2}{l}^{2}+
2412\,k{l}^{3}+676\,{l}^{4}+32\,{k}^{3}+202\,{k}^{2}l\\
\indent +301\,k{l}^{2}+
126\,{l}^{3}+8\,kl+4\,{l}^{2})y_3^6\\
\indent -8\, ( 2\,k+l+1 )  ( 256\,{k}^{5}l+1164\,{k}^{4}{l}^{2
}+1842\,{k}^{3}{l}^{3}+1472\,{k}^{2}{l}^{4}+624\,k{l}^{5}+112\,{l}^{6}
+108\,{k}^{5}+725\,{k}^{4}l+1857\,{k}^{3}{l}^{2}\\
\indent +2072\,{k}^{2}{l}^{3}+
1152\,k{l}^{4}+268\,{l}^{5}+72\,{k}^{4}+400\,{k}^{3}l+800\,{k}^{2}{l}^
{2}+646\,k{l}^{3}+212\,{l}^{4}+12\,{k}^{3}+81\,{k}^{2}l+121\,k{l}^{2}\\
\indent +
63\,{l}^{3}+6\,kl+5\,{l}^{2} ) y_3^5\\
\indent +(800\,{k}^{6}l+9792\,{k}^{5}{l}^{2}+27248\,{k}^{4}{l}^{3}+34472\,{k}^{3
}{l}^{4}+23124\,{k}^{2}{l}^{5}+8000\,k{l}^{6}+1120\,{l}^{7}+640\,{k}^{
6}+9168\,{k}^{5}l \\
\indent +37136\,{k}^{4}{l}^{2}+64046\,{k}^{3}{l}^{3}+55996\,{
k}^{2}{l}^{4}+24520\,k{l}^{5}+4240\,{l}^{6}+1712\,{k}^{5}+13190\,{k}^{
4}l+36747\,{k}^{3}{l}^{2}+45785\,{k}^{2}{l}^{3}\\
\indent +27136\,k{l}^{4}+6124\,
{l}^{5}+1108\,{k}^{4}+6824\,{k}^{3}l+14530\,{k}^{2}{l}^{2}+12968\,k{l}
^{3}+4168\,{l}^{4}+264\,{k}^{3}+1526\,{k}^{2}l+2483\,k{l}^{2}\\
\indent +1321\,{l
}^{3}+20\,{k}^{2}+140\,kl+160\,{l}^{2}+4\,l)y_3^4\\
\indent -4\, ( l+1 )  ( 2\,k+l+1 )  ( 256\,{k}^{4}l
+1164\,{k}^{3}{l}^{2}+1684\,{k}^{2}{l}^{3}+1024\,k{l}^{4}+224\,{l}^{5}
+108\,{k}^{4}+725\,{k}^{3}l+1694\,{k}^{2}{l}^{2}\\
\indent +1560\,k{l}^{3}+480\,{
l}^{4}+72\,{k}^{3}+358\,{k}^{2}l+616\,k{l}^{2}+316\,{l}^{3}+12\,{k}^{2
}+61\,kl+70\,{l}^{2}+4\,l ) y_3^3\\
\indent +4\, ( l+1 ) ^{2} ( 40\,{k}^{5}+482\,{k}^{4}l+1318\,{k}
^{3}{l}^{2}+1406\,{k}^{2}{l}^{3}+656\,k{l}^{4}+112\,{l}^{5}+125\,{k}^{
4}+800\,{k}^{3}l+1593\,{k}^{2}{l}^{2}\\
\indent +1156\,k{l}^{3}+276\,{l}^{4}+79\,
{k}^{3}+410\,{k}^{2}l+584\,k{l}^{2}+224\,{l}^{3}+17\,{k}^{2}+84\,kl+67
\,{l}^{2}+k+6\,l ) y_3^2\\
\indent -8\, ( l+1 ) ^{3} ( 2\,k+l+1 )  ( 2\,l+1+3
\,k )  ( 3\,{k}^{2}+12\,kl+8\,{l}^{2}+k+2\,l ) y_3\\
\indent
+4\, ( l+1 ) ^{4} ( l+k )  ( 2\,l+1+3\,k
 ) ^{2}.$\\
Since the coefficients of the polynomial $t(y_3)$ are positive for even degree terms and negative for odd degree terms, all the real solutions (if exist) of $t(y_3)=0$ are positive.

For all $k,l\in \mathbb{Z}^+$, we have
\begin{eqnarray*}
h(0)&=&2\, \left( 2\,k+1 \right) ^{2} \left( 5\,k+1 \right) ^{2} \left( 4\,k+2\,l+1 \right) >0, \\
h(1) & = & 32\,{k}^{5}+16\,{k}^{4}l-24\,{k}^{3}{l}^{2}-20\,{k}^{2}{l}^{3}-4\,k{l}
^{4}+24\,{k}^{4}+8\,{k}^{3}l \\
&& -18\,{k}^{2}{l}^{2}-12\,k{l}^{3}-2\,{l}^{4
}+4\,{k}^{3}-3\,k{l}^{2}-{l}^{3}\\
& = & \left( 2\,k+1 \right)  \left( 4\,k+2\,l+1 \right)  \left( k-l
 \right)  \left( 2\,k+l \right) ^{2},  \\
h(\infty)& \rightarrow & +\infty.
\end{eqnarray*}
Furthermore assume that $k<l$. Then $h(y_4)=0$ have at least two real solutions, one is between $0$ and $1$, the other is bigger than $1$. Therefore there exist at least two solutions of the homogeneous Einstein equations of the form
 $$\{y_1=y_2={\frac { \left( k+1 \right) y_{{4}}}{l{y_{{4}}}^{2}+3\,k+1}}, y_3=\alpha(y_4), y_4\neq 1, y_5=y_6=1\},$$
where $\alpha(y_4)$ is a rational polynomial of $y_4$ with positive values and the corresponding left-invariant Einstein metrics are non-naturally reductive by Lemma~\ref{naturallyred}. In conclusion, we proved Theorem~\ref{theorem}.

By Theorem~\ref{theorem}, for any $k\leq [\frac{n-1}{3}]$, $\Sp(n)$ admits at least two non-naturally reductive left-invariant Einstein metrics which are $\Ad(\Sp(k)\times\Sp(k)\times\Sp(n-2k))$-invariant. That is, Theorem~\ref{maintheorem} follows.

\section{Acknowledgments}
This work is supported by National Natural Science Foundation of
China (No 11571182).

\end{document}